\documentclass{article}
\usepackage[utf8]{inputenc}
\usepackage[]{authblk}

\title{The Maintenance Location Choice Problem for Railway Rolling Stock}
\author{Jordi Zomer, Nikola Bešinović, Mathijs M. de Weerdt, Rob M.P. Goverde}
\affil{Delft University of Technology, Delft, The Netherlands}
\date{October 2020}

\usepackage{subfiles}

\usepackage[utf8]{inputenc}

\usepackage{heuristica}
\usepackage[heuristica]{newtxmath}
\usepackage{fullpage}

\usepackage{longtable}
\usepackage{booktabs, ctable}

\usepackage{float}
\usepackage{placeins}
\usepackage{array}

\usepackage{array}
\newcolumntype{P}[1]{>{\raggedright\arraybackslash}p{#1}}
\newcolumntype{C}[1]{>{\centering\arraybackslash}p{#1}}

\usepackage[font={small,it}]{caption}

\usepackage{lipsum}

\usepackage{enumitem}
\usepackage{graphicx}
\usepackage{subcaption}
\usepackage{lscape}

\usepackage{amsmath}

\usepackage[official]{eurosym}

\usepackage{pdfpages}

\usepackage{algorithm}
\usepackage{algpseudocode}

\usepackage{hyperref}

\usepackage{natbib}
\usepackage{url}
\begin{document}

\maketitle

%\subfile{"00 Abstract"}
%\subfile{"01 Introduction"} 
%\subfile{"02 Literature review"}
%\subfile{"03 Problem description"}
%\subfile{"04 Model development"}
%\subfile{"05 Results"}
%\subfile{"06 Conclusion"}

%\documentclass["main.tex"]{subfiles}

%\begin{document}

\begin{abstract}
Due to increasing railway use, the capacity at railway yards and maintenance locations is becoming limiting to accommodate existing rolling stock. To reduce capacity issues at maintenance locations during nighttime, railway undertakings consider performing more daytime maintenance, but the choice at which locations personnel needs to be stationed for daytime maintenance is not straightforward. Among other things, it depends on the planned rolling stock circulation and the maintenance activities that need to be performed. This paper presents the Maintenance Location Choice Problem (\textbf{MLCP}) and provides a Mixed Integer Linear Programming model for this problem. The model demonstrates that for a representative rolling stock circulation from the Dutch railways a substantial amount of maintenance activities can be performed during daytime. Also, it is shown that the location choice delivered by the model is robust under various time horizons and rolling stock circulations. Moreover, the running time for optimizing the model is considered acceptable for planning purposes.

\end{abstract}

%\end{document}
%\documentclass["main.tex"]{subfiles}

%\begin{document}
\section{Introduction} \label{sec:introduction} 
In many countries, the demand for rail transport increasing. For instance in The Netherlands, the total number of passenger kilometers increased with over 30\% in the last 20 years. To accommodate increasing travel demand, the utilization of the available rolling stock units may increase and more rolling stock units may be added to the network. 

In order for rail transport to function properly, the rolling stock that operates on the railway network needs to receive maintenance on a regular basis. 
Rolling stock includes all units that move over railway tracks, such as locomotives, passenger wagons and freight wagons. Rolling stock units are fixed compositions that consist of various carriages. The aim of maintenance is to ensure that the rolling stock remains available, safe and comfortable for passengers \citep{dinmohammadi2016risk}. Maintenance activities can be divided into two categories: regular maintenance, corresponding to maintenance activities with higher frequencies (every 1 to 14 days) and shorter duration (1-3 hours), and heavy maintenance, corresponding to maintenance types with lower frequencies (every several months or less) and longer duration (up to several days)~\citep{Andres2015MaintenanceNetworks}. Since regular maintenance and heavy maintenance have such a different nature, the way they are planned may also differ: regular maintenance may be performed whenever a rolling stock unit has a planned standstill; heavy maintenance, however, often requires a rolling stock unit to be taken completely out of regular service. The current work focuses on regular maintenance.

In general, the maximum interval between consecutive maintenance activities is governed by strict rules that are imposed by railway authorities. Maintenance activities are carried out at so-called maintenance locations, which are railway yards with maintenance facilities, spread over the network. For a maintenance location to be operational, it is necessary that personnel is stationed at a maintenance location. The number of personnel stationed at each location is the railway undertaking's decision and determines, together with the maintenance location design, the availability of a maintenance location. An operator can choose to open a facility at certain moments of the day or night by stationing personnel at this location. In particular, a distinction can be made between daytime operations (which means that a location is opened during the day) and nighttime operations (which means that a location is opened during the night). This distinction is clearly visible in The Netherlands where maintenance is usually carried out during nighttime because the rolling stock is mostly in use during daytime. 

%The increasing use of the capacity of the railway network leads to two issues for rolling stock maintenance. First, due to the increased utilization of rolling stock units and the fact that increasingly many rolling stock units operate on the network, the complexity of the scheduling of maintenance activities is increasing. This raises the need for tools that automate the maintenance scheduling process. Second, the use of the capacity of maintenance locations during nighttime is under pressure. As a result, the Dutch railways are considering to perform more maintenance activities during daytime. This raises the question at which locations maintenance teams needs to be stationed to perform daytime maintenance, referred to as the maintenance location choice. 
Due to increasing use of the railway network, the capacity of maintenance locations during nighttime is under pressure and reaching its limits. As a result, the Dutch Railways (Nederlandse Spoorwegen, NS) are considering performing more maintenance activities during daytime. This raises the question at which existing maintenance locations teams need to be stationed to perform daytime maintenance, referred to as the maintenance location choice. 
These maintenance location choice decisions can be made several months before operations, when the rolling stock circulation, which describes all planned rolling stock movements for a given period, for the next period is usually fixed. 
%Note that 
The challenge is to choose maintenance locations such that more maintenance activities can be performed during daytime without modifying the rolling stock circulation.

In planning, commonly three types can be distinguished being strategical (long-term, typically multiple years before operations), tactical (medium-term, from week to one year before operations) and operational (short-term, within several days before operations). When a plan has to be adjusted during operations due to certain disturbances or disruptions, then it is referred to as (real-time) traffic management. In the problem under consideration, it is assumed that the rolling stock circulation is already fixed, and hence it concerns a tactical planning problem. Currently, the choice regarding which locations to open for daytime maintenance is made at a moment in time when the rolling stock circulation has already been determined, justifying the assumption of a fixed rolling stock circulation.

%edited:
This research introduces the Maintenance Location Choice Problem (\textbf{MLCP}), which formulates an optimal maintenance location choice to facilitate timely regular maintenance of all rolling stock, without changing their planned rolling stock circulation. 
%A maintenance schedule and maintenance location choice are considered to be optimal if they achieve the most important goal NS is currently facing: reducing the amount of work that needs to be performed during nighttime. 
The aim of the model is to reduce the amount of work that needs to be performed during nighttime. The model performance is demonstrated on realistic instances provided by the Dutch railways (NS). 

The main contributions of the current paper are the following:
\begin{itemize}
    \item It defines a new mathematical problem formulation for finding  locations for daytime maintenance to reduce nighttime activities without altering the planned rolling stock circulation.
    \item This addresses the location choice problem at a tactical level.
    \item It investigates a new concept of daytime and nighttime maintenance for rolling stock.
    \item It demonstrates the performance of the proposed model on  real-life experiments on the Dutch railway network.
\end{itemize}

Section~\ref{sec:litrev} examines the scientific literature. Section~\ref{sec:problemdescription} introduces the problem and Section~\ref{sec:modeldevelopment} develops a model for solving it. The model's performance is investigated using various scenarios in Section~\ref{sec:results}. Section~\ref{sec:conclusion} gives the key takeaways and identifies some interesting directions for future research.

%\end{document}
%\documentclass["main.tex"]{subfiles}

%\begin{document}
\section{Literature review} \label{sec:litrev}
%Section~\ref{sec:litrev_sched} discusses relevant scientific literature on rolling stock maintenance scheduling. 
This section %\ref{sec:litreview_maintlocation} 
discusses papers on rolling stock maintenance location choice. It covers both railway and aviation papers, given the systematic similarities of the two domains related to vehicle (rolling stock and aircraft) maintenance. 
%At the end of Sections~\ref{sec:litrev_sched} and~\ref{sec:litreview_maintlocation} a selection of corresponding literature from the field of aviation is discussed, which is relevant due to the systematic similarities with rolling stock maintenance scheduling. 
The characteristics of papers are summarised and the differences are mapped out. The existing gaps are defined in the final part of the section.

% In planning, commonly three types can be distinguished being strategical (long-term, typically multiple years before operations), tactical (medium-term, from week to one year before operations) and operational (short-term, within several days before operations). When a plan has to be adjusted during operations due to certain disturbances and disruption, then it is referred as to (real-time) traffic management. 
% The problem addressed in this paper falls under tactical planning.

%\subsection{Maintenance location choice} \label{sec:litreview_maintlocation}
\citet{Tonissen2019MaintenanceUncertainty} aimed at building new maintenance locations in the railway network, which represents a decision problem at the strategic level. Their focus was on determining optimal maintenance locations under unknown/uncertain train lines and fleet size (and mix). They proposed a two-stage stochastic mixed integer programming model, in which the first stage is to open a facility, and in the second stage to minimize the routing cost for the first-stage location decision for each line plan scenario.

\citet{Tonissen2018EconomiesStock} built on \citet{Tonissen2019MaintenanceUncertainty} by including recovery costs of maintenance location decisions, unplanned maintenance, multiple facility sizes and economies of scale (providing that a location twice as big is not twice as expensive). Since, as a result, the second-stage problem becomes NP-hard, an algorithm was provided with the aim to avoid having to solve the second stage for every scenario. 

\citet{Canca2018TheSystems} considered the simultaneous rolling stock allocation to lines and choice for depot locations in a rail-rapid transit context. They proposed a Mixed Integer Linear Programming (MILP) formulation which appears hard to solve. Therefore they proposed a three-step heuristic approach determining first the minimum number of vehicles needed for each line, subsequently the routes of rolling stock on each line, and lastly the circulation of rolling stock on lines over multiple days together with the depot choice.

\citet{Zomer2019VerkenningBehandelen} considered the railway network in The Netherlands and designed a simulation model to estimate the expected effects of carrying out daytime maintenance at given maintenance locations using historical data. Although the use of historical rolling stock data has advantages, since it accurately describes reality (as opposed to rolling stock planning data), it cannot be used for situations in the future for which only planning data is available. Also, although simulation offers the opportunity to investigate effects for various scenarios of daytime maintenance locations, it cannot be used to systematically optimize the maintenance location choice.

In the area of aviation, \citet{Feo1989FlightPlanning} introduced the problem of assigning aircraft to given flights and simultaneously optimizing the number of maintenance facilities. They solved the problem as a minimum cost multi-commodity flow problem. They come up with a two-phase heuristic approach to solve it: firstly, many possible trip patterns for individual aircraft are computed; secondly, the most possible patters are combined to solve the complete problem.

\citet{Gopalan2014TheProblem} is closely related to the work of \citet{Feo1989FlightPlanning}, and assumed routes during the day are given for each aircraft, although these routes are not yet assigned to specific aircraft. Each route needs to be connected to a route on the next day in such a way that the routing passes through a maintenance location with some given periodicity. The objective was to minimize the number of maintenance locations (one of the differences from \citet{Feo1989FlightPlanning}, that considers cost minimization). They proposed four heuristics to solve their problem.

The aforementioned research focused on maintenance location choice. Another area of research in the field of rolling stock maintenance is maintenance scheduling, which is relevant due to the fact that the optimal maintenance location choice typically depends on the schedules for the rolling stock units that need to be maintained. In the literature, maintenance scheduling is often combined with the assignment of rolling stock units (or aircraft) to train trips (or flight legs) such that maintenance constraints are satisfied. Some examples in  railways are \citet{Herr2017F.Scheduling, Andres2015MaintenanceNetworks,Wagenaar2015MaintenanceRailways, VanHovell2019IncreasingLocations} and some examples in aviation are \citet{clarke1997aircraft, gopalan1998aircraft, Sarac2006ARouting, Gopalan2014TheProblem}.
Those problems may be solved using Mixed Integer Linear Programming (MILP) models (e.g. \citet{Herr2017F.Scheduling, Andres2015MaintenanceNetworks, Wagenaar2015MaintenanceRailways}), sometimes extended with solution approaches such as column generation (e.g. \cite{clarke1997aircraft}). Other problems are solved using heuristics (e.g. \citet{gopalan1998aircraft}).
The model developed in this paper also contains important components of maintenance scheduling that incorporate the timely performance of maintenance activities. In particular, time slots are defined during which maintenance actions shall be performed, i.e. maintenance opportunities (introduced in Section \ref{sec:problemdescription}), while the exact start time and end time of actions are not considered. The way in which the current contribution includes these maintenance scheduling components differs from what is typically found in the literature in at least two ways. Firstly, the current contribution incorporates these maintenance scheduling components in conjunction with location choice. Secondly, since the proposed problem is solved on a tactical level, the maintenance scheduling component is addressed assuming a fixed rolling stock circulation.

Table~\ref{tab:lit_overview} summarises the discussed papers. It shows for each paper whether it is written in the aviation (A) or in the railway (R) context, whether it considers maintenance constraints, whether it considers locations choice, whether it considers timing of maintenance for every individual MU, the adopted solution approach and the planning stage to which the problem relates. 
%
%Some more explanation may be necessary on the column indicating maintenance timing of individual rolling stock units. A paper that considers maintenance does not necessarily consider the maintenance timing of maintenance for each rolling stock unit. An example is the work by \citet{clarke1997aircraft}. They do consider maintenance by requiring that each trip path may not exceed some specified length, but the output of their model does not inform the user when and where maintenance on each rolling stock unit could be performed.

Some more explanation may be necessary on the column indicating maintenance timing of individual rolling stock units. Even if a paper considers maintenance (for example by requiring that each trip path shall not exceed some specified length), it need not provide information on where and when each maintenance activity shall be performed. When a paper indeed includes the latter aspect, it is considered to address maintenance timing of individual MUs.

\begin{table}[H]
    \centering
\begin{tabular}{lp{1cm}C{1.2cm}C{1.5cm}C{1.5cm}p{2cm}p{2cm}}
    \FL
      & A/R & \multicolumn{1}{p{1.5cm}}{Maint. considered} & \multicolumn{1}{p{1.5cm}}{Location choice} & \multicolumn{1}{p{1.5cm}}{Indiv. maint. timing}  & Solution approach & Planning stage
     \ML
     \citet{Tonissen2019MaintenanceUncertainty} & R & x & x & & Optimization & Strategic  \\
     \citet{Tonissen2018EconomiesStock} & R & x & x &  & Optimization & Strategic  \\
     \citet{Canca2018TheSystems} & R & & x & & Heuristic & Strategic \\
     \citet{Feo1989FlightPlanning} & A & x & x & & Heuristic & Strategic \\ 
     \citet{Gopalan2014TheProblem} & A  & x & x & & Heuristic & Strategic \\
     \citet{VanHovell2019IncreasingLocations} & R & x & & x & Optimization & Tactical \\ 
     \citet{Zomer2019VerkenningBehandelen} & R & x & x & & Simulation & Tactical \\
     \textit{Current} & R & x & x & x & Optimization & Tactical 
     \LL
\end{tabular}
    \caption{Overview of the aviation (A) and railway (R) maintenance literature discussed in Section~\ref{sec:litrev}.}
    \label{tab:lit_overview}
\end{table}

% The current research is a contribution to the scientific literature since in the following ways.
% \begin{enumerate}
%     \item It offers a simultaneous optimization of rolling stock maintenance scheduling and maintenance location choice. The current paper is unique since it is, to the authors' best knowledge, the only paper that can be classified in both of the last two columns of Table~\ref{tab:lit_overview}.
%     \item It explicitly distinguishes between daytime maintenance and nighttime maintenance. This is relevant for at least the situation in The Netherlands, where nighttime maintenance is standard and recent developments have led the company to investigate daytime maintenance as well. 
%     \item It considers maintenance location choice for The Netherlands, comparable to the work of \citet{Zomer2019VerkenningBehandelen}, but approaches the problem from an optimization perspective instead of from a simulation perspective. 
% \end{enumerate}

This literature review indicates that several aspects have not been addressed in the currently existing literature. 
%First, although problems relating to rolling stock maintenance scheduling and problems relating to maintenance location have been investigated before, both problems have not been solved simultaneously. 
First, although variants of problems relating to rolling stock maintenance location have been investigated before, this problem has not been tackled at a tactical level, i.e. deciding on which maintenance locations to open for given rolling stock circulation. 
Second, no research is available that explicitly distinguishes between daytime maintenance and nighttime maintenance. This distinction is relevant for the situation in the Netherlands, where nighttime maintenance is standard and recent developments have led companies to investigate daytime maintenance as well. Third, although maintenance location choice for the Netherlands has been addressed from a simulation perspective in \citet{Zomer2019VerkenningBehandelen}, no systematic method has been developed yet to find a set of locations for daytime maintenance to minimize the number of nighttime maintenance activities.

%\end{document}
%\documentclass["main.tex"]{subfiles}

%\begin{document}
\section{Problem description} \label{sec:problemdescription}
The problem that is addressed in the current research is defined as the \textit{Maintenance Location Choice Problem} (\textbf{MLCP}).
The input of the \textbf{MLCP} consists of a rolling stock circulation, a set of potential maintenance locations, and a set of maintenance activities that need to be scheduled.
The circulation describes the planned movements of the rolling stock. The maintenance locations are locations where maintenance activities can be carried out. These maintenance activities each have a given duration and constraints on the time between consecutive ones.
The goal of the \textbf{MLCP} is to find an optimal choice of locations used for maintenance, from the set of potential maintenance locations, such that all maintenance activities can be completed in time. Given that rolling stock movements and out of service time in the network can have a direct influence on where each rolling stock unit can be maintained, it is important to consider the rolling stock circulation and some scheduling aspects. The relevant terms and concepts for the \textbf{MLCP} are given in the following.

\begin{figure}[h]
    \centering
    \includegraphics[scale=0.65]{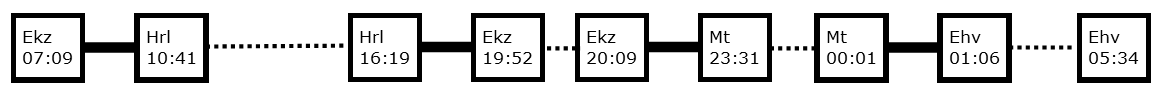}
    \caption{Example of a rolling stock circulation}
    \label{fig:rsc_example}
\end{figure}

\paragraph{Rolling stock circulation} A rolling stock circulation contains for each rolling stock unit a list of trips that a rolling stock unit is scheduled to perform. These trips include an origin and destination position (which are often train stations), and the corresponding planned departure and arrival times. An example of a rolling stock circulation is found in Figure~\ref{fig:rsc_example}. Squares represent arrivals or departures of rolling stock. Inside the squares, the station abbreviation and the departure or arrival time at this station is given. Solid lines represent time intervals where a rolling stock unit is used for a train service. For example, from 07.09 to 10.41, the depicted rolling stock unit is planned to be used for a train service between Ekz (Enkhuizen) and Hrl (Heerlen). Dashed lines represent time intervals where a rolling stock unit is not in service. For example, between 10.41 and 16.19, the depicted rolling stock unit is not in service and will be standing still at Hrl. 

\paragraph{Maintenance opportunities} It is assumed that maintenance can be carried out if and only if a rolling stock unit stands still. These moments are referred to as \textit{maintenance opportuntities} (MOs). Table~\ref{tab:probdes_MOs_example} indicates the MOs corresponding to the rolling stock circulation from Figure~\ref{fig:rsc_example}. All MOs (even short ones) can potentially be used for maintenance activities. When a maintenance activity is assigned to an MO, this implies that the maintenance activity needs to be performed between the start and end time of this MO. The exact scheduled time of the maintenance activity is not determined in the scope of this research.

\begin{table}[h]
    \centering
    \small
\begin{tabular}{llrrr}
    \FL
     MO & day & location & \multicolumn{1}{l}{start time} & \multicolumn{1}{l}{end time}
     \ML
     1 & 1 & Hrl & 10:41 & 16:19 \\
     2 & 1 & Ekz & 19:52 & 20:09 \\
     3 & 1 & Mt & 23:31 & 00:01 \\
     4 & 2 & Ehv & 01:06 & 05:34
     \LL
\end{tabular}
    \caption{Example of the MOs corresponding to the rolling stock circulation example from Figure~\ref{fig:rsc_example}.}\label{tab:probdes_MOs_example}
\end{table}

Two time windows are considered: daytime and nighttime. An MO can be during daytime or nighttime. This division is especially relevant for NS, since it is currently considering a transition from performing maintenance only during nighttime to performing maintenance during both nighttime and daytime \citep{Zomer2019VerkenningBehandelen, VanHovell2019IncreasingLocations}. The time window for daytime maintenance is set at 07.00-19.00, the time window for nighttime maintenance from 19.00-07.00. Some MOs may be partly during daytime and partly during nighttime (this occurs, for example, when an MO lasts from 18.00 to 20.00). For these MOs it is not straightforward whether maintenance would be carried out during the day or during the night (i.e. whether the MO should be considered to be during daytime or during nighttime). For these reasons, the following simplification is applied in the current problem: an MO is marked to be \textit{during daytime} if and only if both its start time and its end time are between 07.00 and 19.00 of the same day; an MO is marked to be \textit{during nighttime} in all other cases. Note that, although this assumption is reasonable in most cases, there are some occasions where it is not realistic: for example, an MO starting at 11.00 (during daytime) and ending at 19.01 (just after the start of the nighttime time window) would be classified to be a nighttime MO whereas maintenance scheduled in it can probably be performed during daytime.

\paragraph{Maintenance activities} It is assumed that a set of regular maintenance types of arbitrary size is known, and that for each maintenance type a fixed duration and a fixed maximum time interval between subsequent maintenance activities is given. For each rolling stock unit, all maintenance types are applied within the respective specified intervals. Specifically, maintenance activities are assigned to maintenance opportunities, which is referred to as the \textit{maintenance scheduling}. It is assumed that maintenance activities are carried out in a subsequent manner and cannot overlap. Furthermore, it is assumed that the number of hours at the start of the planning horizon since the last maintenance activity for each rolling stock unit is given.

\paragraph{Maintenance locations} Maintenance can be carried out at a potential maintenance location. The set of potential maintenance locations is given. For each of these locations, it can be decided whether a location should be opened or not. This decision is referred to as the \textit{location choice}. Each location can be opened during nighttime or during daytime, that is, for each maintenance location, there are four possible outcomes: the location is not used at all, the location is used for daytime maintenance, the location is used for nighttime maintenance, or the location is used for both. It is assumed that each maintenance type can be applied at every location.

%Observe that it is assumed that all maintenance activities can and will be performed at a potential maintenance location, and that there is no subdivision in the types of maintenance that can be performed at specific locations. 

The reachability of maintenance locations is incorporated implicitly since the current research takes the given rolling stock circulation as an input.  Maintenance activities can take place only at those locations where rolling stock units are located according to the rolling stock circulation. 

\paragraph{Planning horizon} The planning horizon in the current research is equal to the planning horizon in the rolling stock circulation. In other words, a maintenance schedule is determined for the entire time horizon of the rolling stock circulation, and as a result, the optimal location choice is valid for the length of this time horizon as well. In The Netherlands, rolling stock circulations are available for periods of eight weeks. This implies that the planning horizon in the current research is also fixed at eight weeks.

\paragraph{Objective} The objective of the \textbf{MLCP} is to assign maintenance activities to MOs and determine for each location whether it is open during daytime and/or nighttime, satisfying the intervals between maintenance activities and the other constraints, in such a way that the number of nighttime maintenance activities is minimized. This goal is relevant in  practical contexts where the capacity of maintenance locations at nighttime is under pressure. As a technical aside, a small penalty applies for any maintenance activity, to avoid the situation that more daytime maintenance activities are planned than necessary.

%\end{document}

%\documentclass["main.tex"]{subfiles}

%\begin{document}

\section{Model development} \label{sec:modeldevelopment}
%The current section introduces a Mixed Integer Linear Programming (MILP) model to solve the \textbf{MLCP}. The use of this type of models enables to find an optimal solution and incorporate various different types of constraints. Moreover, sophisticated solvers exist that can handle problems formulated as a MILP.
The current section provides a Mixed Integer Linear Programming (MILP) model to solve the \textbf{MLCP}. For models formulated as a MILP, sophisticated solvers exist that deliver optimal solutions, and such models have been used for similar problems in the past (e.g. \citet{Tonissen2019MaintenanceUncertainty}). A downside of the use of MILP in general is the fact that the computation time can be relatively high.
The aim therefore is to provide a MILP model that contains not too many decision variables.
Furthermore, given the fact that the \textbf{MLCP} is a problem on the tactical planning level, moderate computation times can be accepted to optimal solutions, as opposed to applications in for instance traffic management (e.g. \citet{Wagenaar2015MaintenanceRailways}), where short computational times are the main goal. 

\subsection{Mathematical notation}
Let $I$ be the set of rolling stock units considered in the current problem and let $i\in I$ be the index used to indicate a specific rolling stock unit.  Let $T > 0$ be the length of the planning horizon in hours. Let $L$ denote the set of potential maintenance locations. 

\paragraph{Maintenance opportunities}
Let $J_i \equiv \{1,..., | J_i |\}$ denote the set of MOs for rolling stock unit $i \in I$. The location of a rolling stock unit $i$ at MO $j \in J_i$ is denoted by $l_{ij} \in L$. The start time of MO $j \in J_i$ is denoted by $s_{ij} \in \mathbb{R}$ and the end time is denoted by $e_{ij} \in \mathbb{R}$, where time is given as the hours that passed since midnight of the first day in the planning horizon. For example, an MO $j$ of a rolling stock $i$ on the fourth day starting at 9:30 and ending at 13:15 would be represented by $s_{ij}=81.5$ and $e_{ij}=85.25$.

Each day is divided into two non-overlapping time windows: the daytime time window and the nighttime time window. Each MO is classified as either a daytime or a nighttime MO. Let $d_{ij}$ indicate whether an MO is classified to be during daytime or during nighttime: let $d_{ij} = 1$ if MO $j \in J_i$ for rolling stock unit $i$ is during daytime and let $d_{ij} = 0$ if a MO $j \in J_i$ for rolling stock unit $i \in I$ is nighttime. Let $\delta^D$ be the hour of the day when the daytime maintenance window starts and let $\delta^N$ be the hour of the day when the nighttime maintenance window starts. Unless stated otherwise, $\delta^D = 7.00$ and $\delta^N = 19.00$. An MO is classified to be during daytime if and only if both its start and end time are during daytime of the same day, i.e.
\begin{align}    
d_{ij} = \begin{cases}
1 & \text{ if } \delta^D \leq e_{ij} \text{ mod } 24 < \delta^N \notag \\
0 & \text{ else} \end{cases}.
\end{align}

\paragraph{Maintenance types}
Let $K$ be the set of maintenance types, $K \equiv \{1, ..., |{K}|\}$. For each maintenance type $k \in K$, let $v_k \in \mathbb{R^+}$ be its duration in hours and let $o_k \in \mathbb{R^+}$ be the maximum interval between two consecutive maintenance activities of maintenance type $k$ in hours.

\paragraph{Maintenance locations}
A potential maintenance location can be opened during daytime (or not) and it can be opened during nighttime (or not). Let $y_l^D \in \{0,1\}$ be a binary variable equal to 1 if location $l \in L$ is available for daytime maintenance and 0 otherwise. Let $y_l^N \in \{0,1\}$ be a binary variable equal to 1 if location $l \in L$ is available for nighttime maintenance and 0 otherwise. The number of potential maintenance locations that can be opened during daytime is restricted. Let $L_{max}^D$ denote the maximum number of potential maintenance locations that can be opened during daytime.

\paragraph{Maintenance schedule}
Maintenance activities are assigned to maintenance opportunities. Let $x_{ijk} \in \{0,1\}$ be a binary variable equal to 1 if maintenance of type $k$ is performed to rolling stock unit $i \in I$ at MO $j \in J_i$, and 0 otherwise. It is required that the total time available at MO $j$ is not exceeded: $\sum_{k \in K} x_{ijk} v_k \leq e_{ij} - s_{ij}$. 

Furthermore, an MO $j$ can only be used if the corresponding location is open at the moment of the MO. Therefore, $d_{ij} = 0$ and $y_{l_{ij}}^N = 0$ implies $x_{ijk} = 0, \forall k \in K$. In the same way, $d_{ij} = 1$ and $y_{l_{ij}}^D = 0$ implies $x_{ijk} = 0, \forall k \in K$.

Moreover, the intervals between successive maintenance activities of the same type should satisfy the given criteria. The interval between two MOs $j, j' \in J_i$, $ j \neq j'$ is measured from the end time of the first MO to the start time of another MO: $s_{ij'} - e_{ij}$. If activity $k\in K$ is scheduled for rolling stock unit $i \in I$ in MO~$j$, then the next maintenance activity should be scheduled such that the interval constraints are satisfied. Let $V_{ijk} \subset J_i$ denote the set of maintenance opportunities for rolling stock unit $i\in I$ that start after the end of MO $j \in J_i$ but earlier than $o_k$ hours after the end of MO $j \in J_i$. This set is $V_{ijk} = \{ p \in J_i : e_{ij} < s_{ip} \leq e_{ij} + o_k \}$ for $j \in J_i$. It is then required that for all $i \in I, j \in J_i$, the following implication holds: 
\begin{align*}
    x_{ijk}=1 \implies \exists j' \in V_{ijk} : x_{ij'k} = 1.
\end{align*}

Observe that a next maintenance activity only needs to be scheduled if maintenance needs to be carried out within the current planning horizon, that is, if ${e_{ij}+o_k \leq T}$.

Let $b_{ik}$ be the number of hours since the last maintenance activity of type $k$ for rolling stock unit $i$ at midnight of the first day. Then let $V_{i0k} = \{ p \in J_i : s_{ip} \leq o_k + b_{ik} \}$.

\subsection{Model formulation}
The \textbf{MLCP} model aims to find $x_{ijk}$ and $y_l$ satisfying the above described constraints that minimize the number of maintenance activities during nighttime. In this model, the decision variables are $x_{ijk}$ ($i \in I, j \in J_i, k \in K$) and $y_l^D$ ($l \in L$). The variables $y_l^N$ ($l \in L$) are considered to be given in the input.

The model can then be formulated as follows:
\begin{align}
    \min \sum_{i \in I} \sum_{j \in J_i} \sum_{k \in K} x_{ijk} (1-d_{ij}) + \varepsilon \sum_{i \in I} \sum_{j \in J_i} \sum_{k \in K} x_{ijk} \label{eq:M11-obj}
\end{align}
subject to
\begin{align}
    1 &\leq \sum_{p \in V_{i0k}} x_{ipk} & \forall i \in I, k \in K \label{eq:M11-c1} \\
    x_{ijk} &\leq \sum_{p \in V_{ijk}} x_{ipk} & \forall i \in I, j \in J_i, k \in K : e_{ij} + o_k \leq T \label{eq:M11-c2}\\
    x_{ijk} &\leq y^D_{l_{ij}} \cdot d_{ij} + y^N_{l_{ij}} \cdot (1-d_{ij}) & \hspace{0.5cm}  \forall i \in I, j \in J_{i}, k \in K \label{eq:M11-c3}\\
    \sum_{k \in K}x_{ijk}v_k &\leq e_{ij} - s_{ij} & \forall i \in I, j \in J_i \label{eq:M11-c4}\\
    \sum_{l \in L} y_l^D & \leq L_{max}^D \label{eq:M11-c5} \\
    x_{ijk} \in \{0,1\}, & \hspace{0.2cm} y_l^D \in \{0,1\} \label{eq:M11-c6}.
   % \sum_{l \in L} y_l^N & \leq L_{max}^N
\end{align}

The objective function (\ref{eq:M11-obj}) minimizes the number of nighttime maintenance activities. The second term penalizes every maintenance activity with an arbitrarily small penalty cost $\varepsilon$ in order to avoid unnecessary maintenance activities being performed. Constraints (\ref{eq:M11-c1}) and (\ref{eq:M11-c2}) enforce that intervals between successive maintenance activities are satisfied. Constraints (\ref{eq:M11-c3}) ensure that maintenance can only be executed at a location that is opened. Constraints (\ref{eq:M11-c4}) take account of the requirement that the duration of maintenance may not exceed the total time of an MO. The number of locations for daytime maintenance is restricted by constraint (\ref{eq:M11-c5}). Constraints (\ref{eq:M11-c6}) ensure that the integer decision variables are also binary.

%\end{document}
%\documentclass["main.tex"]{subfiles}

%\begin{document}
\section{Results} \label{sec:results}
The main goal of the experiments performed in the current section is to investigate the validity of the model by studying the outcomes under various inputs and to make practical recommendations regarding the number of maintenance locations to open and which maintenance locations to open. Also, some numbers regarding computational efficiency are reported.  

%to validate the model (by studying the outcomes) or to (also) say something useful for practice. Also evaluation of computational efficiency?

\subsection{Experimental design} \label{sec:results_expdesign}The proposed \textbf{MLCP} model is tested on real data from the Dutch railway network \citep{NS_BDu}. The data set contains planned rolling stock movements of trips that are operated by NS. A rolling stock circulation covers an 8-week period, which represents a basic planning period, and is issued several weeks before the start of such a period. A data set for such a period is referred to as Basic Day update (in Dutch \textit{BasisDag update}, BDu). This data is an appropriate choice for the current research since it offers a complete view of all rolling stock movements as expected several weeks before operations. In the current research, 10 different BDus have been used. Each of these BDus are individual data sets on which the \textbf{MLCP} model can be run. The period covered by these 10 BDus starts on 10-12-2017 and ends on 1-9-2019. This data set includes all rolling stock units operated by NS, meaning it is of a realistic size and regards a large part of the Dutch rail transport network. The number of rolling stock units in each BDu varies between 820 and 993. Due to the acquisition of new rolling stock units, an increasing trend in the number of rolling stock units operating on the network is visible, leading to the fact that each BDu does not need to contain the same number of rolling stock units.

%Throughout the current results section, 
Two maintenance types are considered ($K = \{1,2\})$: maintenance type 1 has a duration of 0.5 hours ($v_1 = 0.5$) and a maximum interval of 24 hours between successive maintenance activities ($o_1 = 24$), and maintenance type 2 having a duration of 1 hour ($v_2=1$) and a maximum interval of 48 hours between successive maintenance activities ($o_2 = 48$). Rolling stock units are assumed to be as-good-as-new at the start of the planning horizon ($b_{ik} = 0$ for all $i \in I, k \in K$). The start of the daytime time window $\delta^D$ is fixed at $\delta^D = 07.00$ and the start of the nighttime maintenance window $\delta^N$ is fixed at $\delta^N = 19.00$. The technical parameter $\varepsilon$ is fixed at $\varepsilon = 0.001$.

The following parameters are varied to construct the experiments investigated in Section~\ref{sec:results_results}.
\begin{itemize}
    \item We consider 10 BDus $\beta \in \{1, ..., 10 \}$, each representing a period of approximately 8 weeks. Unless stated otherwise, BDu set 10 is used, which is valid in the period 9-6-2019 until 1-9-2019.
    \item The set of rolling stock units $\nu$ is varied. Smaller-sized sets are used to avoid large computation times and thus to allow validating the model by solving multiple instances; larger-sized sets are more realistic which are used to demonstrate a practical performance of the model. In addition, the use of sets of various sizes also enables to investigate the effect of the number of rolling stock units on model outcomes. In most sets, rolling stock units of type VIRM are used since these rolling stock units typically spread out over the entire country visiting many different maintenance locations  and therefore tend to provide most interesting insights on the network level. Specifically, the following sets of rolling stock units are considered:
    \begin{itemize}
        \item Set 1, with the 10 first occurring rolling stock units of type VIRM4.
        \item Set 2, with the 20 first occcurring rolling stock units of type VIRM4.
        \item Set 3, with all rolling stock units of type VIRM4 (in BDu 10: 92)
        \item Set 4, with all rolling stock units of type VIRM4 (in BDu 10: 92) and VIRM6 (in BDu 10: 68)
        \item Set 5, with all rolling stock units used for intercity services (in BDu 10: 360)
    \end{itemize} 
    \item The planning horizon $\tau \in \{7, 21, 42\}$. For meaningful interpretation, $\tau$ should not exceed the length of the period for which the BDu considered is valid (i.e. in general eight weeks). 
    \item The maximum number of daytime maintenance locations  $L^D_{max} \in \{0, 1, 2, 3, 5, 10, 20\}$.
\end{itemize} 

Section~\ref{sec:results_results} presents four experiments to investigate the performance of the \textbf{MLCP} model provided in Section~\ref{sec:modeldevelopment}. 
%exp1
Experiment 1 focuses on the maximum number of locations for daytime maintenance: when more locations can be opened during daytime, supposedly more activities can be performed during nighttime and hence the goal of the \textbf{MLCP} model (minimizing the number of nighttime maintenance activities) may be achieved better.
%exp 2
Experiment 2 focuses on the influence of the planning horizon. If a shorter planning horizon yields a similar location choice, a shorter time horizon may suffice to determine good maintenance locations, which in turn leads to a shorter running time. 
%exp3
Experiment 3 compares 10 different BDus, leading to 10 different scenarios that apply for different time periods, allowing to investigate whether results on the location choice are consistent over multiple time periods.
%exp4
Experiment 4 includes all rolling stock units that serve intercity (long-distance) lines, while the results in the first three experiments considered only the rolling stock types VIRM4 and VIRM6. This allows to investigate the performance and results of the \textbf{MLCP} model in a scenario of realistic size.
More experiments, including experiments elaborating upon the sensitivity of various parameters and the computation time of the \textbf{MLCP} model, can be found in \cite{zomer2020thesis}.
%\textbf{(move inteasibility observation here?
%)
%}

It must be noted that, from the ten available BDus, five resulted in infeasible solutions when applied in combination with rolling stock sets 4 and 5.  This becomes relevant in Experiment 3, where all 10 BDus are considered. These five scenarios are infeasible due to the fact that the corresponding data set contained some rolling stock units for which it is not possible to create a feasible maintenance schedule, for example, no two suitable successive MOs could be found within 24 hours to schedule maintenance activities of type 1. Even if this is the case for only one rolling stock unit, no feasible solution can be found. In reality, this problem could be overcome by manually removing the rolling stock units that caused the infeasible solution and run the model again. However, in the current research, results are reported for the five BDus for which a feasible solution was found only. This is expected to provide sufficient insights as to whether results of one period carry over to multiple periods.

The \textbf{MLCP} is implemented implemented using Python and solved using Gurobi. \citet{GithubMLCP2020} provides the implementation code of the \textbf{MLCP} model used to generate the results in this section. Due to confidentiality, the actual data could not be provided, but a synthetic data set is made available instead.

\subsection{Results} \label{sec:results_results}

\subsubsection{Experiment 1: influence of the maximum number of daytime maintenance locations} \label{sec:results_exp1}

%new:
For varying numbers of the maximum number of daytime maintenance locations, this section presents two measures. Firstly, it presents the \textit{mean daytime share}, which is the percentage of hours of activity performed during daytime, averaged per day, relative to the mean total hours of activity performed per 24 hours, averaged per day. A high mean daytime share indicates that relatively more work is performed during daytime, and this is desired, since it is an indication of a decrease in the pressure on maintenance location capacity during nighttime. However, to check that not only the relative amount of work but also the absolute amount of work during nighttime decreases, it is necessary to verify that the total hours of activity remains more or less constant. Therefore, secondly, this section reports the \textit{mean total hours of activity} per 24 hours, averaged per day, for multiple scenarios.

\paragraph{Mean daytime share}
%To investigate to what extent the goal of the \textbf{MSLCP} is reached, it is relevant to look at the \textit{day share}, which is the percentage of hours of activity performed off average during daytime (relative to the total hours of activity performed on average per day). 
%An aggregate of this statistic is presented in Table~\ref{tab:results_results_exp1_dayshare}.
%Table~\ref{tab:results_results_exp1_dayshare} presents the aggregated results for the day share over multiple scenarios.

Scenarios have been run for all combinations of $L^D_{max}$, $\nu$ and $\tau$, with $L^D_{max} \in \{0, 1, 2, 3, 5, 20 \}$, $\nu \in \{1, 2, 3, 4\}$ and $\tau \in \{ 7, 21, 42\}$ and mean daytime shares for each of these scenarios have been computed. Since the mean daytime shares showed to be relatively invariant for different values of $\nu$ and $\tau$, Table~\ref{tab:results_results_exp1_dayshare} shows \textit{average} mean daytime shares, averaged over all 12 combinations of $\nu$ and $\tau$, for various values of $L^D_{max}$. It shows that the average mean daytime share is increasing with $L^D_{max}$. 

The fact that the mean daytime shares are increasing in $L^D_{max}$ is expected. When more locations are opened for daytime maintenance, more maintenance can be performed during daytime and hence the mean daytime share increases. This is particularly interesting since the goal of NS is to move hours of activity from nighttime to daytime. Apparently, when opening 20 locations, up to on average 42.0\% of the work can be performed during daytime. 

\begin{table}[H]
\small
\centering
\begin{tabular}{lp{2cm}}
\FL
$L^D_{max}$ & average mean daytime share         \ML
0     & \multicolumn{1}{r}{0.0\%} \\
1     & \multicolumn{1}{r}{4.2\%}  \\
2     & \multicolumn{1}{r}{8.9\%}   \\
3     & \multicolumn{1}{r}{13.3\%}  \\
5     & \multicolumn{1}{r}{22.3\%}  \\
20    & \multicolumn{1}{r}{42.0\%}  \LL
\end{tabular}
\caption{The mean daytime share increases for increasing values of $L^D_{max}$. } \label{tab:results_results_exp1_dayshare}
\end{table}

\paragraph{Mean total hours of activity}
%In the previous paragraph, it is indicated that the day share increases for an increasing maximum number of daytime maintenance locations $L^D_{max}$. This is a positive sign, since the objective is to reduce nighttime maintenance activities. 

\begin{figure}[H]
    \centering
    \includegraphics[scale=0.7]{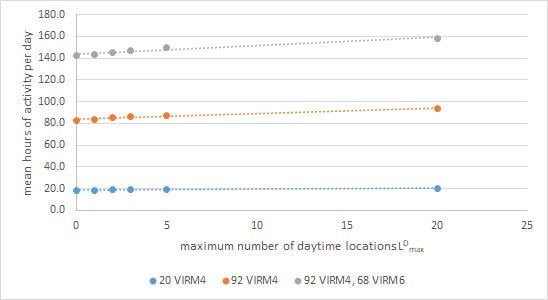}
    \caption{The mean total hours of activity increases when the number of daytime maintenance locations $L^D_{max}$ increases.}
    \label{fig:results_results_exp1_ldmaxhrs}
\end{figure}

%It is also important to investigate the mean total hours of activity. After all, it may be possible that, although the mean day share increases, the total hours of activity increases as well, and this would be undesirable. 
Scenarios have been run for all combinations of $L^D_{max}$, $\nu$ and $\tau$, with $L^D_{max} \in \{0, 1, 2, 3, 5, 20 \}$, $\nu \in \{2, 3, 4\}$ and $\tau \in \{ 7, 21, 42\}$. The set $\nu = 1$ yields similar results but has been excluded from the graph for visualization purposes. The mean total hours of activity for each of these scenarios have been computed. Since the mean total hours of activity showed to be relatively invariant over different values of planning horizon $\tau$, Figure~\ref{fig:results_results_exp1_ldmaxhrs} shows mean total hours of activity per day, averaged over all 3 choices of $\tau$, for various values of $L^D_{max}$ and $\nu$.

Figure~\ref{fig:results_results_exp1_ldmaxhrs} shows the mean total hours of activity for various sets of rolling stock units $\nu \in \{2, 3, 4 \}$. Numbers are averaged over scenarios with different values of planning horizon $\tau \in \{7, 21, 42 \}$.  It shows that indeed the mean total hours of activity per day increases when the maximum number of daytime maintenance locations increases: up to about 15\% when 20 daytime maintenance locations are opened. This implies that, when more possibilities for daytime maintenance arise, this leads to more maintenance in total. From this figure it also becomes clear that this holds for various numbers of rolling stock units, and (evidently) that the mean total hours of activity also increases when more rolling stock units are added to the analysis. The increase of the mean total hours of activity is however convincingly outweighed by the increase of the mean daytime share (cf. Table~\ref{tab:results_results_exp1_dayshare}), and hence the application of daytime maintenance decreases the pressure on the capacity of maintenance locations during nighttime.

\subsubsection{Experiment 2: influence of the planning horizon}
The goal of the second experiment is to investigate the sensitivity of the solutions provided by the \textbf{MLCP} model to various choices of the planning horizon $\tau$. A shorter planning horizon requires less computation time. If a shorter planning horizon yields similar location choice, then the required computation time to determine the optimal location choice could be significantly reduced.

\begin{figure}[H]
    \centering
    \begin{subfigure}[b]{0.48\textwidth}
    \centering
        \includegraphics[scale=0.74]{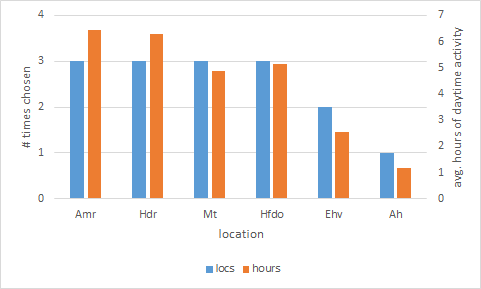}
        \caption{$L^D_{max} = 5$} \label{fig:results_results_exp2_a}
    \end{subfigure}
    \begin{subfigure}[b]{0.48\textwidth}
        \centering
        \includegraphics[scale=0.74]{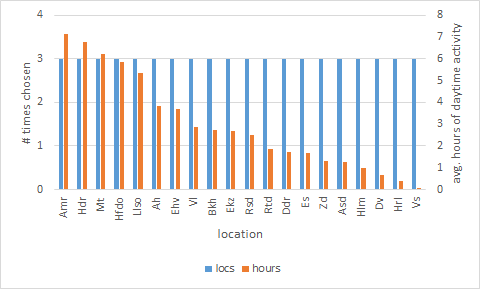}
        \caption{$L^D_{max} = 20$} \label{fig:results_results_exp2_b}
    \end{subfigure}

    \caption{Location choice is consistent over various choices for the planning horizon. Moreover, significant workloads are assigned to various locations. } 
    \label{fig:results_results_exp2}
\end{figure}

Scenarios have been run for all combinations of $L^D_{max}$, $\nu$ and $\tau$, with $L^D_{max} \in \{5, 20 \}$ and $\tau \in \{ 7, 21, 42\}$ with $\nu = 4$ (i.e. all VIRM4 and VIRM6), which means there are three scenarios for each choice of $L^D_{max}$. 

%Figure~\ref{fig:results_results_exp2}, in blue color, shows for each location that was opened in at least one of these three scenarios, in how many of the scenarios this location is opened. This is indicated for $L^D_{max} = 5$ on the left side and for $L^D_{max}=20$ on the right side. 

%new:
Figure~\ref{fig:results_results_exp2} indicates for each location in how many of the investigated scenarios the location was chosen. Moreover, for each location the mean hours of daytime activity assigned to each location is indicated. provides insight in the location consistency and the workload assigned to each location. Blue represents the number of scenarios (out of 3) in which the respective location was open during daytime. Orange represents the average daily number of hours of daytime activity assigned to the respective location.

%The left side of Figure~\ref{fig:results_results_exp2} (blue lines) shows that, when $L^D_{max} = 5$, the locations Amr, Hdr, Hfdo and Mt are chosen in all three scenarios. 
For $L^D_{max} = 5$ (blue lines in Figure~\ref{fig:results_results_exp2_a}), the locations Amr, Hdr, Hfdo and Mt are chosen in all three scenarios. 
This justifies that these four locations are good choices for daytime maintenance in the given subset, but more importantly it shows that the location choice is robust under various time horizons. For $L^D_{max} = 20$ (blue lines in Figure~\ref{fig:results_results_exp2_b}), the location choice is for all three values of planning horizon $\tau$ exactly equal, i.e. each time the same 20 locations are chosen.

Moreover, Figure~\ref{fig:results_results_exp2} shows in orange color the mean hours of \textit{daytime} activity, averaged over all 3 combinations of $\tau$. When $L^D_{max} = 5$ and all rolling stock units are included, the hours of daytime activity is higher than 2 hours per day for these locations that are chosen in all 3 scenarios. In other words, for the locations for which the choice is most consistent, also the highest workload is found.

When comparing Figure~\ref{fig:results_results_exp2_a} and Figure~\ref{fig:results_results_exp2_b}, the locations Amr, Hdr, Mt and Hfdo appear to be four good candidate locations to be opened. They are chosen in all three scenarios, and the workload assigned to these locations is highest. Hence, also when more locations can be opened (20 in this case), the locations found for a smaller maximum number of daytime maintenance locations (5 in this case, see Figure~\ref{fig:results_results_exp2}) still remain good candidates. This is important in practical situations where the number of locations is gradually expanded: the locations found in a stage where only a few locations can be opened, still remain good candidates even when more locations can be opened in a later stage. This is not an obvious finding, for it may have been true that the location choice would change drastically for an increasing value of $L^D_{max}$.

\subsubsection{Experiment 3: consistency over different BDus}

\begin{figure}[H]
    \centering
    \includegraphics{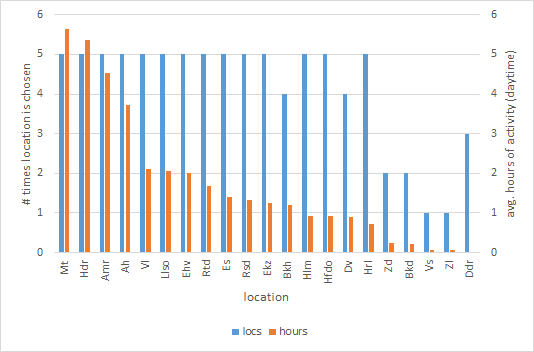}
    \caption{Location choice is consistent over various BDus.}
    \label{fig:results_results_exp3}
\end{figure}

%Scenarios have been run for  $L^D_{max} = 20$, $\nu = 4$ (i.e. all VIRM4 and VIRM6) and $\tau = 42$ for 10 different BDus. From the 10 scenarios run, five resulted in feasible solutions. The other scenarios were infeasible due to the fact that the corresponding data set contained some rolling stock units for which it is not possible to create a feasible maintenance schedule, for example, no two suitable successive MOs could be found within 24 hours to schedule maintenance activities of type 1. Even if this is the case for only one rolling stock unit, no feasible solution can be found. In reality, this problem could be overcome by manually removing the rolling stock units that caused the infeasible solution and run the model again. However, in the current research, results are reported for the five BDus for which a feasible solution was found only. This is expected to provide sufficient insights as to whether results of one period carry over to multiple periods.
Scenarios have been run for  $L^D_{max} = 20$, $\nu = 4$ (i.e. all VIRM4 and VIRM6) and $\tau = 42$. As indicated, out of the 10 availables BDus, only 5 resulted in a feasible solution and hence only these have been considered here.

%(introduce complete figure first. then explain specific results. See Exp 2.)
Figure~\ref{fig:results_results_exp3}, shows for each location that was opened in at least one of scenarios (with a feasible solution), in how many of these five scenarios this location is opened (blue color) and the mean hours of daytime activity at each opened location, averaged over the scenarios in which the location was opened (orange color).

%\textbf{In Figure~\ref{fig:results_results_exp3}, in blue color, it is demonstrated for each location that was opened in at least one of the five scenarios with a feasible solution, in how many of the five scenarios this location is opened.} 
It shows that most of the locations (14 in total) are chosen in all five scenarios (blue lines). 
%for which feasible solutions were obtained. 
This implies that the optimal locations are consistent across different BDus. From a managerial perspective, this is important, since it means that once a location is chosen, it will most likely still be optimal to choose this location in a next period, i.e. for a new rolling stock circulation plan. This means that the applicability of the results, obtained by applying the \textbf{MLCP} model, are in this case not only valid for one BDu, but also carry over to multiple BDus and hence to multiple periods of time. 

In addition, in some cases, a location was chosen in only a limited number of scenarios: for example, Vs and Zl were chosen in only one out of five scenarios. However, in general, these locations are not the locations to which high workloads are assigned. %It may therefore not be of high importance to the ultimate goal of the \textbf{MSLCP} to not open this location.
The influence of choosing not to open such a location may therefore be limited.

%Moreover, Figure~\ref{fig:results_results_exp3} demonstrates in red the mean hours of activity on each of the opened locations, averaged over the five BDus for which a feasible solution was found.
Looking at the mean hours of daytime activity of the opened locations (orange lines), it can be observed that there are four locations with more than three hours of activity and that these locations are consistently chosen over the multiple scenarios. This is positive since it indicates that the choice on whether a location is opened consistently or not correlates with the workload assigned to that location. In other words, if a location is opened consistently throughout multiple scenarios, this usually also means that a relatively high workload is assigned to this location.

\subsubsection{Experiment 4: Performance of scenario of realistic size} \label{sec:results_exp4_comprehensive}
Experiments 1-3 have investigated the performance of the \textbf{MLCP} model under various parameter settings. The largest data set considered so far, included rolling stock units of one train type, VIRM. To assess the workings of the \textbf{MLCP} model in practice, it is necessary to consider a larger and more representative data set. To this end, this fourth experiment considers all rolling stock units used for intercity services combining multiple train types (rolling stock set $\nu = 5$).

Table~\ref{tab:results_results_exp4} gives the most important results for the scenarios with 10 and 20 locations available: it provides the mean hours of daytime activity, the mean daytime share and the running time of the four different scenarios investigated in the current experiment. %from \textbf{the third scenario batch}. 
It shows that, for $L^D_{max}=20$, a mean daytime share of over 30\% can be attained. The mean daytime share appears to be higher for the scenarios with a maximum number of 20 locations opened for daytime maintenance, which is consistent with the findings in Experiment 1 (in Section \ref{sec:results_exp1}).

\begin{table}[H]
\centering
\small
\begin{tabular}{lllrrr}
\FL
$\nu$ & $L^D_{max}$ & mean total hrs. activity & mean daytime share & running time (s) \ML
7          & 10        & 305.9                                        & 22.6\%                                       & 204                                 \\
           & 20        & 310.6                                        & 30.9\%                                     & 201                                 \\
42         & 10        & 328.8                                        & 22.2\%                                      & 11,522                               \\
           & 20        & 334.9                                        & 30.1\%                                     & 7,885          \LL                     
\end{tabular}
\caption{In a realistic case, a substantial amount of work can be performed during daytime.} \label{tab:results_results_exp4}
\end{table}

\begin{figure}[H]
    \centering
    \includegraphics{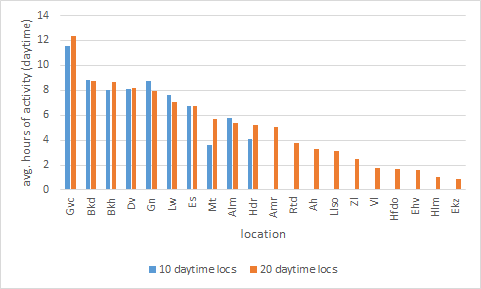}
    \caption{Good candidate locations for daytime maintenance.}
    \label{fig:results_exp7_locuse}
\end{figure}

Figure~\ref{fig:results_exp7_locuse} provides the mean hours of daytime activity per location, for planning horizon $\tau = 42$, set of rolling stock units $\nu = 5$, for $L^D_{max} = 10$ (blue) and $L^D_{max} = 20$ (orange). It shows that, in the situation with a maximum number of locations for daytime maintenance of 10, at least 5 locations have a workload of more than 8 hours per day (Gvc, Bkd, Bkh, Dv and Gn), which can be considered substantial since it is enough to provide work to one maintenance team. Moreover, it shows that the addition of maintenance locations does not seem to reduce the average workload on any of the initial 10 locations. Hence, the initial 10 locations are still good choices, even when daytime maintenance is possible at more locations. 
%new:
Instead, it can be assumed that some additional maintenance work was allocated to these newly opened locations for $L_D^{max}$. 
The added locations, however, are assigned a much lower workload than the initial 10 locations. It may be therefore questionable whether the addition of these locations is worthwhile.

%\end{document}

%\documentclass["main.tex"]{subfiles}

%\begin{document}
\section{Conclusion} \label{sec:conclusion}
%\subsection{Key takeaways}

This paper introduces the Maintenance Location Choice Problem (\textbf{MLCP}) and models it as a mixed integer program, taking as its main input a rolling stock circulation and provides for this rolling stock circulation  an optimal maintenance location choice that minimizes the total number of maintenance activities during nighttime, thereby reducing the capacity pressure during the night. 
For the NS case with all rolling stock units of type VIRM4 and VIRM6, up to 22.3\% of the work can be performed during daytime if five locations are opened, to even 42.0\% if 20 locations are opened during daytime. The location choice is consistent for different lengths of planning horizons and for different input data sets. The four locations with the highest assigned workload (in hours of daytime activity) are Amr, Hdr, Mt and Hfdo. For the largest scenario (including all rolling stock units used for intercity services), a mean daytime share of at least 30.1\% of the activities can be achieved if 20 maintenance locations are opened. In this case, Gvc, Bkd, Bkh en Dv show to be the locations with the highest workloads. Moreover, the model is efficient for planning purposes: for the largest scenario run (with all rolling stock units used for intercity services) for a planning horizon of 42 days, the computation time was 3.2 hours. For operational purposes, however, this running time is not sufficiently efficient. 

%The \textbf{MSLCP} is a contribution to the literature since it simultaneously addresses maintenance scheduling and maintenance location choice. Moreover, its explicit distinction between daytime and nighttime maintenance is an addition to the literature.

Some limitations to the current approach can be identified.
First, the current research uses planned rolling stock data. Although this is common practice in the railway industry and therefore a logical choice, it does not take into account the potential occurrence of disruptions. Due to disruptions, the rolling stock circulation is often altered during operations and as a result, the provided optimal maintenance schedule and optimal maintenance location choice may not be the best anymore at the moment of execution. The current model can provide a basis for more robust decisions on daytime maintenance locations and/or the online adjustment of a maintenance schedule. 
Second, the current research uses a simplified definition of maintenance types. Their durations and intervals may not exactly correspond to practical situations. Moreover, in the current research these are assumed to be equal for each rolling stock unit, whereas in reality these values often differ per rolling stock type. We expect that enriching the formulation with these types would not be very complex and lead to similar results.
Third, the capacity of maintenance locations is not restricted. According to expert judgment at NS, this is acceptable: currently, the primary goal of NS is to reduce the pressure on maintenance locations during nighttime and hence it is not relevant to restrict the capacity of maintenance locations during nighttime. Also, the work assigned to daytime maintenance locations seems to be not excessive and hence it is deemed unnecessary to restrict the capacity of daytime maintenance locations. However, there may be applications in which it is desirable to restrict the capacity of maintenance locations as well, which is an interesting direction for further work.

The \textbf{MLCP} model is intended for use on a tactical level (i.e.\ months before operations), since it assumes that a rolling stock circulation is already given. Hence, it answers the question at which already existing locations personnel needs to be stationed during daytime and/or night time. However, when applied on multiple rolling stock circulations, the \textbf{MLCP} model can also serve more strategic goals (i.e.\ years before operations), regarding for example the question at what positions new locations could best be built.

%\subsection{Future research} \label{sec:conclusion_futureresarch}
Several directions for future research can be recommended.
%First, an interesting next research topic is how to improve the computational performance of the \textbf{MSLCP}. There are several opportunities to improve the computational performance of the \textbf{MSLCP}. Using the structure of the problem, the problem may potentially be decoupled into multiple smaller sub-problems that are much easier to solve in at least three ways. A first opportunity for decoupling may lie in the fact that currently a schedule for all rolling stock units is created simultaneously, while their interaction may be limited. A second opportunity can flow from the fact that the schedule for the several maintenance types is currently created simultaneously for all maintenance types, while the maintenance types possibly do not interact much. A third opportunity may be offered by considering a \textit{rolling horizon} framework, thereby first optimizing a few days in ahead and iteratively adding more days to the optimization. This method would consider only a subset of the decision variables initially and gradually add more decision variables, hence in essence being a column generation approach. Moreover, in the context of the \textbf{CMSLCP}, defining a so-called \textit{warm start} of the solving procedure of the \textbf{MSLCP} model based on its previous outcome may be beneficial.
First, an interesting next research topic is how to improve the computational performance of the \textbf{MLCP}. 
%There are several opportunities to improve the computational performance of the \textbf{MSLCP}. 
%First, using 
Using the structure of the problem, the problem may potentially be decoupled into multiple smaller sub-problems that are much easier to solve, e.g.\ solving rolling stock types independently, or considering a rolling horizon framework, thereby first optimizing a few days in ahead and iteratively adding more days.
%A first opportunity for decoupling may lie in the fact that currently a schedule for all rolling stock units is created simultaneously, while their interaction may be limited. A second opportunity can flow from the fact that the schedule for the several maintenance types is currently created simultaneously for all maintenance types, while the maintenance types possibly do not interact much. A third opportunity may be offered by considering a \textit{rolling horizon} framework, thereby first optimizing a few days in ahead and iteratively adding more days to the optimization. This method would consider only a subset of the decision variables initially and gradually add more decision variables, hence in essence being a column generation approach. Moreover, in the context of the \textbf{CMSLCP}, defining a so-called \textit{warm start} of the solving procedure of the \textbf{MSLCP} model based on its previous outcome may be beneficial.
%Second, the current research can be extended by validating its solutions, which are based on planned data, against realised data. It is interesting to investigate how the solutions provided by the \textbf{MLCP} model may perform in practice. This can potentially be done by computing an \textbf{MLCP} model solution based on planned data from a past period, and investigating its expected performance using realised data from the same period. %This may, for example, be done by constructing a simulation of the historical situation in case a maintenance schedule provided by the \textbf{MSLCP} would have been adopted.
% Third, 
Second, as indicated, some of the input data sets appeared to result in infeasible solutions. Since the model is not capable of allowing (minor) violations of constraints, for some input sets no solution can be found at all. However, in some contexts it may be better to obtain a slightly violated solution than obtaining no solution at all. This could for instance be achieved by accepting and penalising constraint violations such as reducing the time for some maintenance activities by a small amount.
%Fourth, 
Third, the \textbf{MLCP} model can be applied to other domains as well. %It is interesting to extend the models proposed in the current research in order for them to become more versatile. 
Examples of such extensions include the applicability of the research to serve other railway undertakings with potentially different objectives and policies, the investigation of other cost structures and the inclusion of exchange opportunities of rolling stock units so that more rolling stock units can be maintained during daytime.

%\end{document}

\bibliography{References}

\end{document}